\documentclass{article}
\usepackage{amsthm}
\usepackage{amsfonts}

\begin{document}

Tsemo Aristide

College Boreal

1 Yonge Street Toronto, ON

M5E 1E5

tsemo58@yahoo.ca

\bigskip
\bigskip

\centerline{\bf Gerbes, uncertainty and quantization.}

\bigskip
\bigskip

{\it \`A mon p\`ere; tu vivras toujours dans mon coeur.}

\bigskip

\centerline{\bf Abstract.}

\bigskip

The explanation of the photoelectric effect by Einstein and Maxwell's field theory of electromagnetism  have motivated De Broglie to make the hypothesis that matter     exhibits  both waves and particles like-properties. These representations of matter are enlightened by string theory which represents particles with stringlike entities. Mathematically, string theory can be formulated with a gauge theory on loop spaces which is equivalent to the differential geometry of gerbes. In this paper, we show that the descent theory of Giraud and Grothendieck can enable to describe the wave-like properties of the matter with the quantization of a theory of particles. The keypoint is to use the fact that the state space is defined by a projective bundle over the parametrizing manifold which induces a gerbe which represents the geometry obstruction to lift this bundle to a vector bundle. This is equivalent to saying that a phase is determined up to a complex number of module $1$.  When this uncertainty occurs, we cannot locate precisely the position of a particle, and the smallest dimensional quantity that can be described is a $1$ dimensional  manifold; this leads to the concept of wave properties of the matter and string theory.

\bigskip
\bigskip

\centerline{\bf 1. Introduction.}

\bigskip

The structure of matter is a question which has always interested scientists and philosophers.  According to Aristotle, matter was infinitely divisible,   while Democritus proposed that matter has a building block, an idea which has leaded to the concept of atom.

The point of view of Democritus has been adopted by modern scientists until the experiences of diffraction of Fresnel and Young   confirmed the wave theory of Huygens about the light. But, in 1905, Einstein has introduced the notion of photon to explain the photoelectric effect of the light. These facts have motivated De Broglie to assert that every particle has wave-like and particle-like properties. This duality can be used to explain the uncertainty principle of Heisenberg  which asserts that the position and the momentum of a particle cannot be simultaneously determined since it is difficult to locate a wave on a string.

In this paper, we use descent theory defined by Giraud and Grothendieck to formulate  a geometric uncertainty principle which enable to interpret the results of experiences which exhibit the wave  nature of the matter. The observations of the matter depend of how precise is defined the state space and henceforth the corresponding phase space. In this paper, we study quantum systems whose group of symmetries is the circle group $U(1)$, and which are in an adiabatic evolution.  The phase space is parametrized by a manifold on which is defined the Hamiltonian. If the Hamiltonian is defined by an operator on an Hilbert bundle over the parametrized manifold, the  Berry phase  of the system is described by the holonomy of a connection. The $2D$-uncertainty principle asserts that we cannot determine the state space of  the quantum system precisely; in mathematical terms, it is defined by a gerbe on the parametrizing manifold. This  is equivalent to saying that,  locally we can only determine the isomorphism class of the quantum system. In this case, we can only compute the Berry phase up to an error: the $2D$-Berry phase which is the holonomy of the gerbe mentioned above along a surface. 

As Aristotle, we do not believe in the existence of a fundamental particle which plays the role of the building block of the matter. The development of the technology enables to make very sophisticated experiences which show the existence of new particles often smaller than the previous known. But we believe that a theory which aims to describe (part of) the matter as to start somewhere by choosing a fundamental structure. Since the technology that we use to describe the matter is not perfect, we have to take in account uncertainty in the theory that we use. In so doing, the objects that our theories can describe maybe different from the initial  fundamental structure and represent a set of them. In fact, the transgression morphism identifies an $U(1)$-gerbe defined on a manifold $M$ to a line bundle defined on its loop space $LM$. With this identification, the $2D$-Berry phase is the holonomy of a connection defined on $LM$, instead of determining the position of  a fundamental particle, we can only determine the position of a set of these particles which is a string. To explain the theory of photons,  Einstein in [4] says the following:  the larger the energy density and the wavelength of radiation the more suitable the theoretical basis that we used; but for small wavelengths and low radiation densities the basis fails completely. When the energy density is big, it can be measured with higher precision the uncertainty can be neglected, if it is small, we have to take in account uncertainty of measurements.

\bigskip

{\bf Notation.}

\medskip

Let $M$ be a manifold, and $C$ a presheaf of categories defined on $M$. Let $(U_i)_{i\in I}$ be an open covering. For every subset $\{i_1,...,i_n\}$ of $I$, we will often denote the intersection of $U_{i_1},...,U_{i_n}$ by $U_{i_1..i_n}$. Let $e_{i_1...i_n}$ be an object of $C(U_{i_1...i_n})$, for every subset $\{j_1,...,j_m\}$ of $I$, we will denote by $e_{i_1..i_n}^{j_1..j_m}$ the restriction of $e_{i_1...i_n}$ to $U_{i_1..i_nj_1..j_m}$. If $f_{i_1..i_n}:e_{i_1...i_n}\rightarrow e'_{i_1..i_n}$ is a morphism between objects of $C(U_{i_1..i_n})$, we will denote by $f_{i_1..i_n}^{j_1..j_m}:e_{i_1..i_n}^{j_1..j_m}\rightarrow {e'}_{i_1..i_n}^{j_1..j_m}$ its restriction to $U_{i_1..i_nj_1..j_m}$.

 \bigskip

\centerline{\bf 2. $1D$-quantum systems and the  Berry phase.}

\bigskip

In this section we are going to define quantum systems and  review the definition of the Berry phase of a $1D$-quantum system in an adiabatic evolution. The notions of gerbe and connection defined  on gerbe used in this paper are defined in the appendix.

 \medskip
 
 {\bf Definitions 2.1.}

 A $2D$-quantum space is defined by:
 
  A complex Hilbert space  ${\cal H}$;
  
  The state space, which  is the projective space $P({\cal H})$ associated to ${\cal H}$.
  
  A $P({\cal H})$-bundle $p_M:P\rightarrow M$.

 Let $U({\cal H})$ be the unitary group of ${\cal H}$ the extension of groups:
  
  $$
  1\longrightarrow U(1)\rightarrow U({\cal H})\rightarrow PGl({\cal H})\rightarrow 1
  $$
  
 defines a gerbe $C_{\cal H}$ on $M$ bounded by  $U_M(1)$, the sheaf of $U(1)$ valued functions defined on $M$ such that for every open subset $U$ of $M$, $C_{\cal H}(U)$ is the category whose objects are ${\cal H}$-bundles over $U$ whose associated projective bundle is the restriction of $p_M$ to $U$.
 
 We supposed defined on every object $e_U$  of $C_{\cal H}(U)$ an Hamiltonian operator $H_{e_U}$, that  is an operator of the ${\cal H}$-bundle $e_U$. We denote by $D_{\cal H}(U)$ the subcategory of $C_{\cal H}(U)$ which has the same class of objects than $C_{\cal H}(U)$ and such that $Hom_{D_{\cal H}}(e_U,e'_U)$ is the subset of $Hom_{C_{\cal H}}(e_U,e'_U)$ whose elements commute with the respective Hamiltonian.
 
 We assume that the correspondence $U\rightarrow D_{\cal H}(U)$ is a gerbe defined on $M$ bounded by the sheaf of $U(1)$-functions defined on $M$.
 
Suppose that the gerbe $D_{\cal H}$ is trivial, a $1D$-quantum system is a global section $(P,H,{\cal H})$ of $D_{\cal H}$.  Where $P$ is the corresponding global section of $C_{\cal H}$ and $H$ the Hamiltonian of $P$.

  \medskip
  
  {\bf Definition 2.2.}

  We denote by $Aut(P,{\cal H},H)$ the group of automorphisms of the $1D$-quantum system $(P,{\cal H},H)$. Let $Diff(M)$ be the group of diffeomorphisms of $M$. There exists a morphism $\pi:Aut(P,{\cal H},H)\rightarrow Diff(M)$ defined by the following commutative diagram. 
  
$$
\matrix{P&{\buildrel{g}\over{\longrightarrow}} & P\cr p_M\downarrow && \downarrow p_M\cr M &
{\buildrel{\pi(g)}\over{\longrightarrow}}& M}
$$

  \medskip

 The Schroedinger equation on $(P,{\cal H},H)$ is defined by:

$$
i{\partial\over{\partial t}}\Psi(t,x)=H\Psi(t,x) 
$$

Here $\Psi$ the wave function, we suppose that $|\Psi|=1$.

\medskip

Suppose that $\Psi(0,x_0)$ is in an eigenvector of $H(x_0)$. The quantum adiabatic theorem says that if a quantum system is perturbed slowly, it remains in its eigenstate if there is a gap between the eigenvalue and the rest of the spectrum of the Hamiltonian.

To apply this result, we suppose that for every $x \in M$, $H_x$ has an eigenspace $E(x)$ of dimension $1$. The initial condition $\Psi(t_0,x_0)\in E(x_0)$. We also assume that the family $E(x), x\in M$ is a subbundle $p_E:E\rightarrow M$ of $p_M:P\rightarrow M$. This bundle inherits an Hermitian metric.We will call the bundle $p_E:E\rightarrow M$ the energy level bundle.

 We suppose that there exists  an atlas $(U_j)_{j\in J}$ of $M$ such that for each $j\in J$,   there exists a section $\phi_j$ of $E$ such that  $\langle\phi_j,\phi_j\rangle=1$. 
 
 The quantum adiabatic theorem enables us to  look local solutions of the form:
  
$$
\Psi_j =e^{i\theta_j}\phi_j.
$$

The Schroedinger equation gives:

$$
i{\partial\over{\partial t}}\Psi_j(t,x)=-d\theta_je^{i\theta_j}\phi_j+ie^{i\theta_j}d\phi_j=E(x)e^{i\theta_j}\phi_j
$$

By making the Hermitian product of the last equation with $\phi_j$, we deduce that:

$$
-d\theta_j+i\langle d\phi_j,\phi_j\rangle=E_j(x)
$$
since $\langle \phi_j,\phi_j\rangle=1$.

We write:

$$
\alpha_j=-i\langle d\phi_j,\phi_j\rangle.
$$

\medskip

On $U_j\cap U_k$, we can write:

$$
\phi_k^j=e^{ia_{jk}}\phi^k_j
$$
and the family $(e^{ia_{jk}})_{j,k\in J}$ is family of transition functions of the bundle $E$.

We   have: 
 
 $$
 {\alpha_k}^j=-i\langle d(e^{ia_{jk}}\phi^k_j),e^{ia_{jk}}\phi^k_j\rangle
 $$
 
 $$
 =-i\langle ida_{jk}e^{ia_{jk}}\phi^k_j,e^{ia_{jk}}\phi^k_j\rangle-i\langle e^{ia_{jk}}d\phi^k_j,e^{ia_{jk}}\phi_j^k\rangle
 $$
 
 $$
 =\alpha^k_j+da_{jk}.
 $$

 This is equivalent to saying that the family of $1$-forms  $\alpha_j$ define the  Koszul derivative of a connection $\nabla_{\cal H}$ defined on the energy level bundle since the transition functions of $p_E$ are $(e^{ia_{jk}})$.  We denote by $\omega_{\cal H}$ the Ehresmann connection associated to $\nabla_{\cal H}$ defined on the $U(1)$-principal bundle associated to $E$. The connection $\omega_{\cal H}$ will be called an admissible connection The holonomy of the connection $\omega_{\cal H}$  is the Berry phase.

\bigskip

{\bf 3. $2D$-quantum systems.}

\bigskip

In this section, we do not suppose that the gerbe $C_{\cal H}$ is trivial. This is equivalent to assume that we can only determine the position of a point up to an element of $U(1)$. As a consequence of this uncertainty, we cannot determine precisely the Berry phase around a loop. We are going to see that the gerbe $D_{\cal H}$ is endowed with a connection whose holonomy can be interpreted as a Berry phase on the loop space of $M$.

\medskip

Let $D_{\cal H}$ be a $2D$-quantum system, we suppose that the rank of the energy level bundle $E_{e_U}$ of every object $e_U$ of $D_{\cal H}(U)$ is $1$. Let $(U_j)_{j\in J}$ be a good covering of $M$, and $e_j$ an object of $D_{\cal H}(U_j)$. We denote by $p_j:e_j\rightarrow U_j$ the bundle map. We denote by $\omega_j$ a fundamental connection of the energy level bundle $E_j$ of $e_j$. Let $g_{jk}:e_k^j\rightarrow e_j^k$ a morphism of $D_{\cal H}(U_{jk})$. The classifying cocycle of $D_{\cal H}$ is defined by $c_{jkl}=g_{lj}^kg_{jk}^lg_{kl}^j$.

We define:

$$
\omega_{jk}=\omega_j^k-g^*_{kj}\omega_k^j.
$$

There exists a form $\alpha_{jk}$ defined on $U_{jk}$ such that $\omega_{jk}{=p_j^k}^*\alpha_{jk}$.

We have:

$$
\alpha_{kl}^j-\alpha_{jl}^k+\alpha_{jk}^l=i dlog(c^{-1}_{jkl}).
$$

The family of forms $\alpha_{jk}$ is a connction defined on $D_{\cal H}$. It defines a connective structure on $D_{\cal H}$ see the proof of proposition 5.3.2  in [2]. We call that connective structure the fundamental connective structure.

 \medskip

 {\bf Schrodinger equation and $2D$-quantum systems.}

 \medskip
 
 Let $D_{\cal H}$ be a $2D$-quantum system defined on the manifold $M$. An interesting problem is to determine the solutions of the Schrodinger equation. In particular, we are interested to study the Berry phase. But it is not possible to determine precisely the position of a particle. We are going to show that it is possible to determine the error induced by the uncertainty principle. This error is the holonomy of a gerbe around a surface.

 \medskip
 
Let $c:I=S^1\rightarrow M$ be a differentiable curve. We want to compute the Berry phase along $c$. The pullback $c^*D_{\cal H}$ is a trivial gerbe, but may have more than one global sections. Let $e_c$ and $e'_c$ be global sections of $c^*D_{\cal H}$, they are endowed with an admissible connection, but the holonomy of each of these connections depends of the global object. 

Let $c:[0,1]\times S^1$ be a differentable map. The pullback  $c^*D_{\cal H}$ by $c$ is a trivial gerbe.The following result is proved in the appendix:

$$
\int_{[0,1]\times S^1}Hol(c^*C^M_{\cal H})=hol(c_0)hol(c_1)^{-1}
$$

where $c_i:[0,1]\rightarrow M$ is defined by $c_i(t)=c(i,t), i=0,1$, $hol(c_i)$ is the holonomy of an admissible connection on a global object of $c_i^*D_{\cal H}$, and $Hol(c^*C^M_{\cal H})$ is the holonomy of the pullback of the fundamental connective structure by $c$.

\medskip

The fundamental object studied in physics is the smallest  object that the theory can describe. Here it is a  string since for a $2D$-quantum system, we cannot determine the Berry phase of a particle on a loop, but we can measure the evolution of  the Berry phase on a string. 

Let $LM$ be the loop space of $M$, the transgression morphism induces an isomorphism between $H^2(LM,{\mathbb Z})$ and $H^3(M,\mathbb{Z})$. In fact there exists a canonical correspondence between $U(1)$-gerbes endowed with a connective structure and $U(1)$-line bundles defined on $LM$ endowed with a connection. With this identification, the Berry phase defined on a loop on $LM$  corresponds to the holonomy of a surface of a gerbe defined on $M$.

\bigskip

{\bf Appendix Gerbes.} 

\bigskip

{\bf Definition A1.}
 
 Let $M$ be a differentiable manifold  a sheaf of categories is a correspondence $C$ defined on $Ouv(M)$ the set of open subsets of $M$, which associates to every open subset $U$ of $M$ a category $C(U)$ such that:

 - For  every open subsets $U,V$ of $M$ such that $U\subset V$, there exists a restriction functor $r_{U,V}:C(V)\rightarrow C(U)$ such that:
 
 $$
 r_{U,V}\circ r_{V,W}=r_{U,W}
 $$
 
 Gluing condition for objects.
 
 Let $U$ be an open subset of $M$ and $(U_i)_{i\in I}$ an open covering of $U$. Let $e_i$ be an object of $C(U_i)$, we write $e_i^j=r_{U_i\cap U_j,U_i}(e_i)$. Let $u_{ij}:e_j^i\rightarrow e_i^j$ be a family of morphisms such that:
 
 $$
 u_{ij}^ku_{jk}^i=u_{ik}^j
 $$
 there exists an object $e_U$ of $C(U)$ such that $r_{U_i,U}(e_U)=e_i
 $.
 
 Gluing condition for morphisms.
 
 Let $e,e'$ be objects of $C(U)$ write $e_i=r_{U_i,U}(e)$ and ${e'}_i=r_{U_i,U}(e')$. Let $f_i:e_i\rightarrow {e'}_i$ be morphisms of $C(U_i)$ such that:
 
 $$
 f_j^i=f_i^j.
 $$
 There exists a morphism $f:e\rightarrow e'$ such that $r_{U_i,U}(f)=f_i$.
 
  The sheaf of categories $C$ is called a gerbe if and only if:
  
  The category $C(U)$ is a groupoid;
  
  For every $x\in M$ there exists a an open neighbourhood $U_x$ of $x$ such that $C(U_x)$ is not empty;
  
  Let $e_U,e'_U$ be objects of $C(U)$, for every element $x$ of $U$, there exists an open subset $U_x\subset U$ containing $x$ such that $r_{U_x,U}(e_U)$ is isomorphic to $r_{U_x,U}(e'_U)$. 
  \bigskip
  
 {\bf Definition A2.}
  
  The gerbe $C$ is bounded by the sheaf $L$ defined on $M$ if and only if for every object $e_U$ of $C(U)$, there exists an isomorphism $l_{e_U}:Aut(e_U)\rightarrow L(U)$ between $Aut(e_U)$ the group of endomorphisms of $e_U$ and $L(U)$ which commutes with  restrictions and with the morphisms between objects.
  
\bigskip

{\bf The classifying cocycle.}

\medskip

We denote by $U_M(1)$  the sheaf of $U(1)$ differentiable functions defined on $M$. We are going to associate to a gerbe bounded by $U_M(1)$ a classifying $2$-Cech cocycle.

Let $(U_i)_{i\in I}$ be a good open covering of $M$ such that for every $i\in I$, $C(U_i)$ is a non empty connected groupoid, and $e_i$ an object of $C(U_i)$, there exists a morphism:

$$g_{ij}:e_j^i\rightarrow e_i^j
$$

The morphism $c_{ijk}=g_{ki}^jg_{ij}^kg_{jk}^i$ is an automorphism of $e_k^{ij}$, we can identify it with an element of $U_M(1)(U_{ijk})$ by using $l_{e_k^{ij}}$. We have:

$$
c_{ikl}^jg_{lk}^{ij}c_{ijk}^lg_{kl}^{ij}=c_{ijl}^kc_{jkl}^i
$$

The morphism $c_{ikl}^j, c_{ijl}^k,c_{jkl}^i$ are automorphisms of $e_l^{ijk}$ which can be identified to an element of $U_M(1)(U_{ijkl})$. The morphism $g_{lk}^{ij}c_{ijk}^lg_{kl}^{ij}$ is also a morphism of $e_l^{ijk}$ which is identified with the same element of $U_M(1)(U_{ijkl})$ than the morphism $c_{ijk}^l$ of $e_k^{ijl}$ (see definition A2). We deduce that:

$$
c_{ijk}^l(c_{ijl}^k)^{-1}c_{ikl}^l(c_{jkl})^{-1}=1.
$$

 The theorem of Giraud implies that there exists a bijection between the set of equivalence gerbes bounded by $U_M(1)$ and the Cech-cohomology group $H^2_{Cech}(M,U_M(1))$.

\bigskip

\bigskip

{\bf Connections on  gerbes.}

\medskip

{\bf Definition A3.}

 Let $M$ be a manifold, and $C$ an $U_M(1)$-gerbe. We suppose that an object $e_U$ of $C(U)$ is a fibre bundle $e_U\rightarrow U$ endowed with a free action of $U(1)$. Remark that we do not assume that $e_U$ is a principal $U(1)$-bundle defined on $U$. A morphism $f:e_U\rightarrow e'_U$ is a morphism of fibre bundles which commutes with the action of $U(1)$. We can identify the band of $C$ with the sheaf of $U(1)$-valued functions defined on $M$.
 
 Let $X$ be an element of the $Lie(U(1))\simeq \mathbb{R}$, there exists a vector field $X_{e_U}$ defined on $e_U$ such that $X_{e_U}(x)={d\over{dt}}_{t=0}exp(tX).x$.
 
   A connection $\omega$ defined on $C$ is a correspondence which associates to every object $p_{e_U}:e_U\rightarrow U$ of $C(U)$ a
    $1$-form $\omega_{e_U}$, which takes its values in $Lie(U(1))$ such that:
   
   for every $x\in e_U$ and $\omega_{e_U}(X_{e_U})=X$,
   
   For every $g\in U(1), g^*\omega=\omega$.
   
   Let $f:e_U\rightarrow e'_U$ a morphism of $C$, there exists a $1$-form $\alpha_{e_U,e'_U}$ defined on $U$ such that $\omega_{e_U}-f^*\omega_{e'_U}=p_{e_U}^*\alpha_{e_U,e'_U}.$ 
   
  Remark that for every automorphism $g$ of $e_U$, we have $g^*\omega_{e_U}=\omega_{e_U}-idlog(g)$.

 \bigskip

Let $(U_j)_{j\in J}$ a good open covering of $M$  and $e_j$ an object of $C(U_j)$. Let $g_{jk}:e_k^j\rightarrow e_j^k$
 a morphism.
 
 We define: $
 \omega_{jk}=\omega_j^k-g_{kj}^*\omega_k^j.
 $
 
 We have:
 
 $$
 \omega_{kl}^j-\omega_{jl}^k+\omega_{jk}^l=g_{kj}^*(\omega^{jl}_k-g_{lk}^*\omega_l^{jk})-(\omega_j^{kl}-g_{lj}^*\omega_l^{jk})+\omega_j^{kl}-g_{kj}^*\omega_k^{jl}.
 $$
 
 $$
 =g_{lj}^*(\omega_l^{jk}-(g_{lk}^jg_{kj}^lg_{jl}^k)^*\omega_l^{jk})=idlog(c_{jkl}^{-1}).
 $$
 
 There exists a $1$-form $\alpha_{ij}$ defined on $U_{ij}$ such that $\omega_{ij}=(p_j^i)^*\alpha_{ij}$ where $p_i^j:e_i^j\rightarrow U_{ij}$ is the bundle morphism. We also have:
 
 $$
 \alpha_{kl}^j-\alpha_{jl}^k+\alpha_{jk}^l=idlog(c_{jkl}^{-1}).
 $$
 
 The family of forms $\alpha_{jk}$ is called a connection form defined on the gerbe $C$.

We deduce that the family $d\alpha_{jk}$ defines a cohomology class in $H^1_{Cech}(M,\Omega^2(M))$. Since this group vanishes, there exists a family of $2$-forms $F_i$ defined on $U_i$ such that:

$$
F_k^j-F_j^k=d\alpha_{jk}.
$$

\medskip

{\bf Definition A3.}

{\it The curvature $R_C$ of the gerbe is the $3$-form whose restriction to $U_j$ is $dF_j$.}

\medskip

 Suppose now that $M$ is a $2$-dimensional manifold, since $H^2(M,U(1))=1$, we deduce that the gerbe $C$ is trivial. There exists a family of functions $h_{jk}:U_j\cap U_k\rightarrow U(1)$ such that:
 
 $$
 h_{kl}^jh_{lj}^kh_{jk}^l=c_{jkl}.
 $$
 
 This implies that:
 
 $\alpha_{jk}+idlog(h_{jk})$ is a Cech cocycle. There exists a family of $1$-form $m_j$  defined on $U_j$ such that:
 
 $$
 m^j_k-m^k_j=\alpha_{jk}+id(logh_{jk}).
 $$
 
 The restrictions of the $2$-forms $F_j-dm_j$ and $F_k-dm_k$ coincide on $U_j\cap U_k$. We denote by $Hol(M)$ the $2$-form defined on $M$ whose restriction to $U_j$ is $F_j-dm_j$. The holonomy of the gerbe is:
 
 $$
 exp(i\int_MHol(M)).
 $$
 
  \bigskip

Let $f:N\rightarrow M$ be a differentiable map, we can define the pullback $f^*C$ of $C$. For every open subset $U$ of $M$, an object of $f^*C(f^{-1}(U))$ is the pullback of an object of $C(U)$ by $f$. 

Suppose that the gerbe $C$ is endowed with the connection $\omega$, the gerbe $f^*C$ is endowed with the connection $f^*\omega$   such that for every object $e_U$ of $U$ endowed with the connection $\omega_{e_U}$, $f^*e_U$ is endowed with the connection $f^*\omega_{e_U}$.
 
 We do not assume  that the dimension of $M$ is $2$, let $S$ be a $2$-dimensional manifold and $f:S\rightarrow M$ a differentiable map. We denote by $f^*C$ the pullback of $C$ by $f$. The holonomy of $C$ around $S$ is:
 
 $$
 exp(i\int_SHol(f^*C)).
 $$
 
 \bigskip
 
 {\bf Theorem. A3.}
 
 {\it Let $c:[0,1]\times S^1\rightarrow M$ be a differentiable map, then:
 
 $$
 exp(i\int_{[0,1]\times S^1}Hol(c^*C))=Hol(c_0)Hol(c_1)^{-1}
 $$
 where $hol(c_i), i=0,1$ is the holonomy of a global object of $c^*_iC$ around the loop $c_i:I\rightarrow S^1$  defined by $c_i(t)=c(i,t)$.}
 
\medskip

{\bf Proof.}

We can use a global object $e$ of $c^*C$ to define the holonomy form $Hol(c^*C)$. Since $H^2([0,1]\times S^1,\mathbb{Z})$ is trivial, $e$ is trivial. We deduce that there exists an open covering $(U_i)_{i\in I}$ of $M$, a family of maps $(d_i:U_i\rightarrow U(1))_{i\in I})$  such that the coordinate changes of $e$ are defined by $d_j^i(d_i^j)^{-1}$.

We can use a connection on $e$ defined locally by $(\omega_i)_{i\in I}$ to define the $1$-Cech chain of $1$-forms which defines the connection $c^*C$. Let $\alpha_{ij}$ be the two form defined on $U_i\cap U_j$ whose pullback by $e_i^j\rightarrow U_i\cap U_j$ is we have $\omega_i^j-g_{ji}^*\omega_i^j$. Since $d\alpha_{ij}=0$,  we can suppose the $0$-Cech cocycle of $2$-forms $F_i$ associated to the connection  of $c^*C$ defined by the family $\alpha_{ij}$ is zero. The holonomy cocycle is then $d(-\alpha_i-idlog d_i)$.

 The Stoke's theorem implies  that:

$$
\int_{[0,1]\times S^1}Hol(c^*C)=\int_{[0,1]\times[0,1]}h^*Hol(c^*C)=\int_{h^*c_0}Hol(c_0^*e)-\int_{h^*c_1}Hol(c_1^*e).
$$

This implies that:

 $$
 exp(i\int_{[0,1]\times S^1}Hol(c^*C))=Hol(c_0)Hol(c_1)^{-1}.
 $$

\bigskip

\centerline {\bf References.}

\bigskip

[1]. Berry Michael. {\it Quantals phase factors accompanying adiabatic changes.} Proceedings of the Royal Society A. 392 (1802) 45-57.

\smallskip

[2]. Brylinski J.L; {\it Loop spaces characteristic classes and geometric quantization.} 1993

\smallskip

[3]. De Broglie Louis.{\it  Recherches sur la th\'eorie des quantas} Ann. de Physique 10 (3) 22 1925

\smallskip

[4]. Einstein, A. {\it Concerning an heuristic point of view toward the emission and transformation of light.} American Journal of Physics 33.5 (1965): 367.

\smallskip

[5]. Giraud, J. {\it Cohomologie non abélienne de degré 2.} PhD diss., Impr. administrative centrale, 1966.

\end{document}